\documentclass[a4paper,notitlepage,11pt]{article}%
\usepackage{amssymb}
\usepackage{graphicx}
\usepackage{amsmath}
\usepackage{amsfonts}%
\providecommand{\U}[1]{\protect\rule{.1in}{.1in}}
\numberwithin{equation}{section}
\numberwithin{theorem}{section}
\ifx\pdfoutput\relax\let\pdfoutput=\undefined\fi
\newcount\msipdfoutput
\ifx\pdfoutput\undefined\else
\ifcase\pdfoutput\else
\ifx\paperwidth\undefined\else
\ifdim\paperheight=0pt\relax\else\pdfpageheight\paperheight\fi
\ifdim\paperwidth=0pt\relax\else\pdfpagewidth\paperwidth\fi
\fi\fi\fi
\begin{document}

\title{Existence of strictly positive solutions for sublinear elliptic problems in
bounded domains \thanks{2000 \textit{Mathematics Subject Clasification}.
35J25, 35J61, 35B09, 35J65.}  \thanks{Partially supported by Secyt-UNC. } }
\author{T. Godoy, U. Kaufmann \thanks{\textit{E-mail addresses. }godoy@mate.uncor.edu
(T. Godoy), kaufmann@mate.uncor.edu (U. Kaufmann, Corresponding Author).}
\and \noindent\\{\small FaMAF, Universidad Nacional de C\'{o}rdoba, (5000) C\'{o}rdoba,
Argentina}}
\maketitle

\begin{abstract}
Let $\Omega$ be a smooth bounded domain in $\mathbb{R}^{N}$ and let $m$ be a
possibly discontinuous and unbounded function that changes sign in $\Omega$.
Let $f:\left[  0,\infty\right)  \rightarrow\left[  0,\infty\right)  $ be a
nondecreasing continuous function such that $k_{1}\xi^{p}\leq f\left(
\xi\right)  \leq k_{2}\xi^{p}$ for all $\xi\geq0$ and some $k_{1},k_{2}>0$ and
$p\in\left(  0,1\right)  $. We study existence and nonexistence of strictly
positive solutions for nonlinear elliptic problems of the form $-\Delta
u=m\left(  x\right)  f\left(  u\right)  $ in $\Omega$, $u=0$ on $\partial
\Omega$.

\end{abstract}

\section{Introduction}

Let $\Omega\subset\mathbb{R}^{N}$, $N\geq1$, be a $C^{1,1}$ bounded domain.
Our aim in this paper is to consider the question of existence of solutions
for nonlinear problems of the form%
\begin{equation}
\left\{
\begin{array}
[c]{ll}%
-\Delta u=mf\left(  u\right)  & \text{in }\Omega\\
u>0 & \text{in }\Omega\\
u=0 & \text{on }\partial\Omega,
\end{array}
\right.  \label{prob}%
\end{equation}

where $m:\Omega\rightarrow\mathbb{R}$ is a function that changes sign in
$\Omega$ and $f:\left[  0,\infty\right)  \rightarrow\left[  0,\infty\right)  $
is a continuous function satisfying

H1. $f$ is nondecreasing, and there exist $k_{1},k_{2}>0$ and $p\in\left(
0,1\right)  $ such that $k_{1}\xi^{p}\leq f\left(  \xi\right)  \leq k_{2}%
\xi^{p}$ for all $\xi\geq0$.

As pointed out in \cite{nodea}, the existence of strictly positive solutions
for sublinear problems with indefinite nonlinearities as (\ref{prob}) raises
many interesting questions and is intriguing even in the one-dimensional case
for various reasons. One of them is that the existence of (nontrivial)
nonnegative solutions does not guarantee the strict positivity of such
solutions (in contrast for example to superlinear problems, where they even
belong to the interior of the positive cone). In fact, there are situations in
which there exist nonnegative solutions which actually vanish in a subset of
$\Omega$ (see e.g. \cite{bandle}). Another one is for instance that several
\textit{non-comparable} sufficient conditions on $m$ can be established for
the existence of solutions for (\ref{prob}) in the one-dimensional case under
some evenness assumptions on $m$ (see \cite{nodea}, Section 2), and these
solutions may not be in the interior of the positive cone.

The present work is a natural continuation of the research started in
\cite{nodea}, where $m$ was considered (when $N>1$) to be radially symmetric.
Let us note that the nonlinearity studied there was $f\left(  \xi\right)
=\xi^{p}$. One of the most important differences between $\xi^{p}$ and the
nonlinearities treated in this paper is that here (\ref{prob}) is no longer
homogeneous in $m$ (i.e. (\ref{prob}) may admit a solution but $km$ may not
($k>0$ constant), and viceversa), and the homogeneity was crucial in the every
existence proofs given in \cite{nodea}.

We shall primarily rely on the well-known sub- and supersolution method in the
presence of weak sub and supersolutions (see e.g. \cite{du}, Theorem 4.9). One
of the reasons is that the existence of supersolutions represent no
difficulty, see Remark 2.3 below. In order to supply (strictly positive)
subsolutions, we shall divide the domain in parts and construct subsolutions
in each of them, and later check that they can be joined appropriately to get
a subsolution in the entire domain. This last fact depends on obtaining
estimates for the normal derivatives of these subsolutions on the boundaries
of the subdomains. In \cite{nodea} these bounds could be computed rather
explicitly making use of the radial symmetry of $m$ (and the fact that
$\Omega$ was a ball) but in the present situation those computations cannot be
done any more. Let us mention that here the key tool will be an estimate due
to Morel and Oswald, see Lemma 2.1 below.

In Theorem 3.1 we shall state a sufficient condition on $m$ for the existence
of solutions of (\ref{prob}), while in Theorem 3.2 we shall provide a
\textquotedblleft local\textquotedblright\ necessary condition and a
\textquotedblleft global\textquotedblright\ one in Corollary 3.3 under an
additional assumption on $m$. We observe that this last condition is of
similar type to the one in Theorem 3.1. In order to relate these results to
others already existing, we mention that two necessary conditions were proved
for some particular radial functions in \cite{nodea}, Theorem 3.4 (see also
Remark 3.5 there), and as far as we know there are no other results (other
than the obvious condition $m^{+}\not \equiv 0$ implied by the maximum
principle). Concerning the matter of sufficient conditions, the only theorem
we found in the literature, apart from the ones proved in \cite{nodea} for $m$
radial, is that there exists a solution for (\ref{prob}) provided that the
solution of the linear problem $-\Delta\phi=m$ in $\Omega$, $\phi=0$ on
$\partial\Omega$, satisfies $\phi>0$ in $\Omega$ (see \cite{jesusultimo},
Theorem 4.4, or \cite{hand}, Theorem 10.6). As a matter of fact, this even
holds for linear second order elliptic operators with nonnegative zero order
coefficient. We note however that the aforementioned condition is far from
being necessary in the sense that there are examples of (\ref{prob}) having a
solution but with the corresponding $\phi$ satisfying $\phi<0$ in $\Omega$
(cf. \cite{nodea}). Let us finally mention that for $m$ smooth an $p\in\left(
0,1\right)  $ it is known that the problem $-\Delta u=mu^{p}$ in $\Omega$,
$u=0$ on $\partial\Omega$ admits a (nontrivial) nonnegative solution if and
only if $m\left(  x_{0}\right)  >0$ for some $x_{0}\in\Omega$ (see e.g.
\cite{bandle} or \cite{publi}).

We conclude this introduction with some few words on the case of a general
second order elliptic operator. We believe that at least some of the results
presented here should still be true when $-\Delta$ is replaced by such
differential operators. In fact, one can verify that except the use of Lemma
2.1, the proof of Theorem 3.1 can be carried out exactly as it is done here
(with the obvious changes) in the case of a general operator. Hence, if a
similar version of the aforementioned lemma holds for these operators (which a
priori it is not clear since the proof makes use of the mean value properties
for superharmonic functions), then an analogue of Theorem 3.1 can be proved in
this case.

\textit{Acknowledgments}. The authors are pleased to thank the referee for
her-his careful and detailed reading of the paper.

\section{Preliminaries}

The following estimate appeared first in an unpublished work by Morel and
Oswald (\cite{morel}), and a nice proof can be found in the paper of Brezis
and Cabr\'{e}, \cite{cabre}, Lemma 3.2.

$\qquad$

\textbf{Lemma 2.1. }\textit{Let }$h\in L^{r}\left(  \Omega\right)  $\textit{,
}$r>N$\textit{, and let }$u$\textit{ be the solution of }%
\begin{equation}
\left\{
\begin{array}
[c]{ll}%
-\Delta u=h & \text{\textit{in} }\Omega\\
u=0 & \text{\textit{on }}\partial\Omega.
\end{array}
\right.  \label{h}%
\end{equation}
\textit{Then there exists some }$c=c\left(  \Omega\right)  >0$\textit{ such
that}
\begin{gather*}
u\left(  x\right)  \geq c\delta_{\Omega}\left(  x\right)  \int_{\Omega}%
h\delta_{\Omega}\qquad\text{\textit{for all }}x\in\Omega,\\
\text{\textit{where }}\delta_{\Omega}\left(  x\right)  :=dist\left(
x,\partial\Omega\right)  .
\end{gather*}

$\qquad$

The next result is also known (see e.g. Theorem 3.4 in \cite{uriel2}). We
present a brief sketch of the proof for the sake of completeness. Let us note
that the following proof is much simpler than the one given in \cite{uriel2}.
We set
\[
P^{\circ}\overset{.}{=}\text{interior of the positive cone of }C^{1,\alpha
}\left(  \overline{\Omega}\right)  \text{, }\alpha\in\left(  0,1\right)
\text{.}%
\]

$\qquad$

\textbf{Lemma 2.2. }\textit{Let }$m\in L^{r}\left(  \Omega\right)  $\textit{
with }$r>N$\textit{ and such that }$0\not \equiv m\geq0$, \textit{and let }%
$f$\textit{ satisfying H1. Then there exists }$v\in W^{2,r}\left(
\Omega\right)  \cap P^{\circ}$\textit{ solution of }%
\begin{equation}
\left\{
\begin{array}
[c]{ll}%
-\Delta v=mf\left(  v\right)  & \text{in }\Omega\\
v>0 & \text{in }\Omega\\
v=0 & \text{on }\partial\Omega.
\end{array}
\right.  \label{z}%
\end{equation}
\textit{Proof}. Let $\phi>0$ be the solution of $-\Delta\phi=m$ in $\Omega$
and $\phi=0$ on $\partial\Omega$. Then using the second inequality in H1 one
can verify that for every $k>0$ large enough it holds that $k\left(
\phi+1\right)  $ is a supersolution of (\ref{z}). On the other side, let
$\varphi>0$ with $\left\Vert \varphi\right\Vert _{\infty}=1$ satisfying
\[
\left\{
\begin{array}
[c]{ll}%
-\Delta\varphi=\lambda_{1}\left(  m,\Omega\right)  m\varphi & \text{in }%
\Omega\\
\varphi=0 & \text{on }\partial\Omega,
\end{array}
\right.
\]
where $\lambda_{1}\left(  m,\Omega\right)  $ denotes the (unique) positive
principal eigenvalue for $m$. It is easy to check employing the first
inequality in H1 that $\varepsilon\varphi$ is a subsolution of (\ref{z}) for
all $\varepsilon>0$ sufficiently small, and the lemma follows. $\blacksquare$

$\qquad$

\textbf{Remark 2.3. }Let us mention that the construction of the supersolution
made in the first part of the above proof still works if $m$ changes sign in
$\Omega$, taking there $\phi$ as the solution of $-\Delta\phi=m^{+}$ in
$\Omega$ and $\phi=0$ on $\partial\Omega$. (where as usual we write
$m=m^{+}-m^{-}$ with $m^{+}=\max\left(  m,0\right)  $ and $m^{-}=\max\left(
-m,0\right)  $). Furthermore, this is also true for a strongly uniformly
elliptic differential operator with nonnegative zero order coefficient.
$\blacksquare$

\section{Main results}

\textbf{Theorem 3.1. }\textit{Let }$\Omega_{0}$ \textit{be a }$C^{1,1}$
\textit{domain with }$\overline{\Omega}_{0}\subset\Omega$\textit{, and
let\ }$m\in L^{r}\left(  \Omega\right)  $\textit{ with }$r>N$ \textit{and
}$0\not \equiv m\geq0$\textit{\ in }$\Omega_{0}$\textit{. Let }$k_{1}$,
$k_{2}$ \textit{be given by H1. There exist some }$C_{0},C_{1}>0$
\textit{depending only on }$\Omega$ \textit{and }$\Omega_{0}$\textit{ such
that if}%
\[
\left\Vert m^{-}\right\Vert _{L^{r}\left(  \Omega-\overline{\Omega}%
_{0}\right)  }\leq\frac{k_{1}C_{0}}{k_{2}C_{1}^{1-p}}\int_{\Omega_{0}}%
m\delta_{\Omega_{0}}^{p+1}%
\]
\textit{then }(\ref{prob}) \textit{has a solution }$u\in W^{2,r}\left(
\Omega\right)  $.

\textit{Proof}. Let $\Omega-\overline{\Omega}_{0}:=\Omega_{1}$. For $M>0$, we
start constructing some $0\leq w\in W^{2,r}\left(  \Omega_{1}\right)  $
solution of
\begin{equation}
\left\{
\begin{array}
[c]{ll}%
-\Delta w=-m^{-}f\left(  w\right)   & \text{in }\Omega_{1}\\%
\begin{array}
[c]{l}%
w=0\\
w=M
\end{array}
&
\begin{array}
[c]{l}%
\text{on }\partial\Omega\\
\text{on }\partial\Omega_{0}.
\end{array}
\end{array}
\right.  \label{w}%
\end{equation}
Let us first note that since by H1 $f\left(  0\right)  =0$, it holds that
$\underline{w}:=0$ is a subsolution of (\ref{w}), and also since $f$ is
nonnegative we have that $\overline{w}:=M$ is a supersolution of (\ref{w}). It
follows from Theorem 4.9 in \cite{du} that there exists some $w$ weak solution
of (\ref{w}) satisfying $0\leq w\leq M$. Furthermore, by standard arguments we
may conclude that $w\in W^{2,r}\left(  \Omega_{1}\right)  $ (indeed, it is
enough to note that if $z\in W^{2,r}\left(  \Omega_{1}\right)  $ is the unique
solution of the problem $-\Delta z=-m^{-}f\left(  w\right)  $ in $\Omega_{1}$,
$z=0$ on $\partial\Omega$ and $z=M$ on $\partial\Omega_{0}$, then the maximum
principle implies that $z=w$). 

We claim now that there exists some $C>0$ depending only on $\Omega_{1}$ such
that if $M:=\left[  Ck_{2}\left\Vert m^{-}\right\Vert _{L^{r}}\right]
^{1/\left(  1-p\right)  }$ then $w>0$ in $\Omega_{1}$ ($k_{2}$ given by H1).
To confirm this, let $\theta,\psi\in W^{2,r}\left(  \Omega_{1}\right)  $ be
the unique solutions of%
\[
\left\{
\begin{array}
[c]{ll}%
\Delta\theta=0 & \text{in }\Omega_{1}\\%
\begin{array}
[c]{l}%
\theta=0\\
\theta=1
\end{array}
&
\begin{array}
[c]{l}%
\text{on }\partial\Omega\\
\text{on }\partial\Omega_{0},
\end{array}
\end{array}
\right.  \qquad\quad\left\{
\begin{array}
[c]{cc}%
-\Delta\psi=m^{-} & \text{in }\Omega_{1}\\
\psi=0 & \text{on }\partial\Omega_{1}.
\end{array}
\right.
\]
From the Sobolev imbedding theorems and the $W^{2,r}$-theory for elliptic
equations (e.g. \cite{grisvard}, Theorem 2.4.2.5) we derive that
\[
\left\vert \psi\right\vert \leq\left\Vert \nabla\psi\right\Vert _{L^{\infty}%
}\delta_{\Omega_{1}}\leq\left\Vert \psi\right\Vert _{C^{1}}\delta_{\Omega_{1}%
}\leq c_{0}\left\Vert \psi\right\Vert _{W^{2,r}}\delta_{\Omega_{1}}\leq
c_{1}\left\Vert m^{-}\right\Vert _{L^{r}}\delta_{\Omega_{1}}%
\]
for some $c_{1}=c_{1}\left(  \Omega_{1}\right)  >0$, and we also have that
$\theta>c_{2}\delta_{\Omega_{1}}$ in $\Omega_{1}$ for some $c_{2}=c_{2}\left(
\Omega_{1}\right)  >0$.

On the other hand, since $w\leq M$, recalling H1 we get that in $\Omega_{1}$%

\[
-\Delta\left(  M\theta-k_{2}M^{p}\psi\right)  =-m^{-}k_{2}M^{p}\leq-m^{-}%
k_{2}w^{p}\leq-m^{-}f\left(  w\right)  =-\Delta w
\]
and so
\[
w\geq M\theta-k_{2}M^{p}\psi>\left(  c_{2}M-c_{1}k_{2}M^{p}\left\Vert
m^{-}\right\Vert _{L^{r}}\right)  \delta_{\Omega_{1}}\qquad\text{in }%
\Omega_{1}%
\]
and the claim is proved. We fix for rest of the proof $M$ as in the
aforementioned claim.

Let $\nu$ denote the outward unit normal to $\partial\Omega_{0}$. Let us
observe now that
\begin{gather}
\left\vert \frac{\partial w}{\partial\nu}\right\vert \leq\left\Vert
w\right\Vert _{C^{1}}\leq c_{0}\left\Vert w\right\Vert _{W^{2,r}}\leq
c_{1}\left(  M+\left\Vert m^{-}\right\Vert _{L^{r}}\left\Vert f\left(
w\right)  \right\Vert _{L^{\infty}}\right)  \leq\label{deri}\\
c_{1}\left(  M+k_{2}M^{p}\left\Vert m^{-}\right\Vert _{L^{r}}\right)
\leq2c_{1}\left[  \max\left\{  1,C\right\}  k_{2}\left\Vert m^{-}\right\Vert
_{L^{r}}\right]  ^{1/\left(  1-p\right)  }:=\nonumber\\
c_{3}\left[  c_{4}k_{2}\left\Vert m^{-}\right\Vert _{L^{r}}\right]
^{1/\left(  1-p\right)  },\nonumber
\end{gather}
with $c_{3}$ and $c_{4}$ depending only on $\Omega_{1}$. 

On the other side, let $v>0$ be the solution of (\ref{z}) with $\Omega_{0}$ in
place of $\Omega$. Taking into account H1 and Lemma 2.1, there exists
$c_{5}=c_{5}\left(  \Omega_{0}\right)  >0$ such that $v\geq c_{5}k_{1}%
\delta_{\Omega_{0}}\int_{\Omega_{0}}mv^{p}\delta_{\Omega_{0}}$ and so raising
this inequality to the power $p$, multiplying by $m\delta_{\Omega_{0}}$ and
integrating over $\Omega_{0}$ we obtain $\left(  \int_{\Omega_{0}}mv^{p}%
\delta_{\Omega_{0}}\right)  ^{1-p}\geq\left(  c_{5}k_{1}\right)  ^{p}%
\int_{\Omega_{0}}m\delta_{\Omega_{0}}^{1+p}$ and hence%
\[
v\geq\left[  c_{5}k_{1}\int_{\Omega_{0}}m\delta_{\Omega_{0}}^{1+p}\right]
^{1/\left(  1-p\right)  }\delta_{\Omega_{0}}\text{.}%
\]
Define now $u:=M+v$. Then $\partial u/\partial\nu\leq-\left[  c_{5}k_{1}%
\int_{\Omega_{0}}m\delta_{\Omega_{0}}^{p+1}\right]  ^{1/\left(  1-p\right)  }$
and $u=w$ on $\partial\Omega_{0}$. Hence, if we set $\omega:=u$ in
$\overline{\Omega}_{0}$ and $\omega:=w$ in $\overline{\Omega}-\Omega_{0}$ it
follows applying the divergence theorem (as stated e.g. in \cite{cuesta}, p.
742) that $\omega$ is a weak subsolution of (\ref{prob}) if $\partial
u/\partial\nu\leq\partial w/\partial\nu$. Recalling (\ref{deri}) this occurs
if
\[
c_{3}^{1-p}c_{4}k_{2}\left\Vert m^{-}\right\Vert _{L^{r}}\leq c_{5}k_{1}%
\int_{\Omega_{0}}m\delta_{\Omega_{0}}^{p+1}%
\]
and thus, taking into account Remark 2.3, this ends the proof. $\blacksquare$

\qquad

We denote with $B_{R}\left(  x_{0}\right)  $ the open ball in $\mathbb{R}^{N}$
centered at $x_{0}$ with radius $R$, and we write $\left(  -\Delta\right)
^{-1}:L^{r}\left(  \Omega\right)  \rightarrow L^{\infty}\left(  \Omega\right)
$ for the solution operator of (\ref{h}). We also set
\begin{equation}
C_{N,p}:=\frac{\left(  1-p\right)  ^{2}}{2\left(  N\left(  1-p\right)
+2p\right)  }. \label{cnp}%
\end{equation}

\qquad

\textbf{Theorem 3.2. }\textit{Let }$m\in L^{r}\left(  \Omega\right)  $\textit{
with }$r>N$\textit{, let }$C_{N,p}$\textit{ be given by }(\ref{cnp})
\textit{and let }$k_{1}$, $k_{2}$\textit{ be given by H1. If there exists a
solution }$u\in C\left(  \overline{\Omega}\right)  $ of (\ref{prob}),
\textit{then }%
\begin{gather}
\frac{C_{N,p}}{\left\Vert \left(  -\Delta\right)  ^{-1}\right\Vert }%
\sup_{B_{R}\left(  x_{0}\right)  \in\mathfrak{B}}\left[  m_{R}R^{2}\right]
<\frac{k_{2}}{k_{1}}\left\Vert m^{+}\right\Vert _{L^{r}\left(  \Omega\right)
}\text{,}\qquad\text{\textit{where}}\label{nec}\\
\mathfrak{B}:=\left\{  B_{R}\left(  x_{0}\right)  \subset\Omega:m\leq0\text{
\textit{in }}B_{R}\left(  x_{0}\right)  \right\}  \text{,}\qquad m_{R}%
:=\inf_{B_{R}\left(  x_{0}\right)  }m^{-}\text{.}\nonumber
\end{gather}

\textit{Proof}. We proceed by contradiction. If (\ref{nec}) does not hold,
then there exists some $B_{R}\left(  x_{0}\right)  \in\mathfrak{B}$ such that
\begin{equation}
\frac{C_{N,p}m_{R}R^{2}}{\left\Vert \left(  -\Delta\right)  ^{-1}\right\Vert
}\geq\frac{k_{2}}{k_{1}}\left\Vert m^{+}\right\Vert _{L^{r}\left(
\Omega\right)  }.\label{acsur}%
\end{equation}
Let $\beta:=1/\left(  1-p\right)  $, and for $x\in\overline{B}_{R}\left(
x_{0}\right)  $ define
\[
w\left(  x\right)  :=\left[  k_{1}C_{N,p}m_{R}\left\vert x-x_{0}\right\vert
^{2}\right]  ^{\beta}\text{.}%
\]
After some computations one can verify that $\Delta w\leq k_{1}m^{-}w^{p}$ in
$B_{R}\left(  x_{0}\right)  $. Let $u$ be a solution of (\ref{prob}). In
particular, it holds that $\Delta u\geq k_{1}m^{-}u^{p}$ in $B_{R}\left(
x_{0}\right)  $. Also, taking into account H1, from (\ref{prob}) we deduce
that
\begin{equation}
\left\Vert u\right\Vert _{L^{\infty}\left(  \Omega\right)  }\leq\left[
k_{2}\left\Vert \left(  -\Delta\right)  ^{-1}\right\Vert \left\Vert
m^{+}\right\Vert _{L^{r}\left(  \Omega\right)  }\right]  ^{\beta}.\label{cita}%
\end{equation}
Moreover, if $x\in\partial B_{R}\left(  x_{0}\right)  $, employing
(\ref{acsur}) and (\ref{cita}) we derive that%
\[
w\left(  x\right)  =\left(  k_{1}C_{N,p}m_{R}R^{2}\right)  ^{\beta}\geq\left[
k_{2}\left\Vert \left(  -\Delta\right)  ^{-1}\right\Vert \left\Vert
m^{+}\right\Vert _{L^{r}\left(  \Omega\right)  }\right]  ^{\beta}%
\geq\left\Vert u\right\Vert _{L^{\infty}\left(  \Omega\right)  }\geq u\left(
x\right)  .
\]
It follows by the comparison principle that $w\geq u$ in $B_{R}\left(
x_{0}\right)  $, but $w\left(  x_{0}\right)  =0$, contradicting the fact that
$u>0$ in $\Omega$. $\blacksquare$

\qquad

\textbf{Corollary 3.3. }\textit{Let }$\Omega_{1}\subset\Omega$ \textit{be a
convex domain and let} $m\in L^{r}\left(  \Omega\right)  $\textit{ with }$r>N$
\textit{and such that in }$\Omega_{1}$ $m$ \textit{is convex and }$m\leq0$.
\textit{If there exists a solution }$u\in C\left(  \overline{\Omega}\right)
$\textit{ of }(\ref{prob}),\textit{ then }%
\begin{equation}
\frac{4C_{N,p}}{27\left\vert \Omega_{1}\right\vert \left\Vert \left(
-\Delta\right)  ^{-1}\right\Vert }\int_{\Omega_{1}}m^{-}\delta_{\Omega_{1}%
}^{2}<\frac{k_{2}}{k_{1}}\left\Vert m^{+}\right\Vert _{L^{r}\left(
\Omega-\overline{\Omega}_{1}\right)  }\text{.} \label{coro}%
\end{equation}
\textit{Proof}. Let $\alpha:=2/3$ and let $x_{1}\in\Omega_{1}$. We set
$R_{1}:=\alpha\delta_{\Omega_{1}}\left(  x_{1}\right)  $ and let $y\in
B_{R_{1}}\left(  x_{1}\right)  $. Observe that $z_{y}\left(  t\right)
:=x_{1}+t\left(  y-x_{1}\right)  \in\Omega_{1}$ for every $t\in\left[
0,1/\alpha\right]  $ since $\left\vert z_{y}\left(  t\right)  -x_{1}%
\right\vert <\delta_{\Omega_{1}}\left(  x_{1}\right)  $. Define $M\left(
t\right)  :=m^{-}\left(  z_{y}\left(  t\right)  \right)  $. Then $M\left(
t\right)  $ is concave in $\left[  0,1/\alpha\right]  $ and hence%
\[
m^{-}\left(  y\right)  =M\left(  1\right)  \geq\alpha M\left(  1/\alpha
\right)  +\left(  1-\alpha\right)  M\left(  0\right)  \geq\left(
1-\alpha\right)  M\left(  0\right)  =m^{-}\left(  x_{1}\right)  /3.
\]
It follows that $\inf_{B_{R_{1}}\left(  x_{1}\right)  }m^{-}\geq m^{-}\left(
x_{1}\right)  /3$. Now, if (\ref{prob}) possesses a solution\textit{ }$u\in
C\left(  \overline{\Omega}\right)  $, by Theorem 3.2 we obtain that%
\begin{gather*}
\frac{k_{2}}{k_{1}}\left\Vert m^{+}\right\Vert _{L^{r}\left(  \Omega\right)
}>\frac{C_{N,p}}{\left\Vert \left(  -\Delta\right)  ^{-1}\right\Vert }%
\sup_{B_{R}\left(  x_{0}\right)  \in\mathfrak{B}}\left[  m_{R}R^{2}\right]
\geq\\
\frac{C_{N,p}}{3\left\Vert \left(  -\Delta\right)  ^{-1}\right\Vert }%
m^{-}\left(  x_{1}\right)  R_{1}^{2}=\frac{4C_{N,p}}{27\left\Vert \left(
-\Delta\right)  ^{-1}\right\Vert }m^{-}\left(  x_{1}\right)  \delta
_{\Omega_{1}}^{2}\left(  x_{1}\right)
\end{gather*}
for every $x_{1}\in\Omega_{1}$. Integrating this inequality in $\Omega_{1}$
with respect to $x_{1}$ gives (\ref{coro}) and thus the corollary is proved.
$\blacksquare$

\qquad

\textbf{Remark 3.4. }We observe that $C_{N,p}\rightarrow0$ when $p\rightarrow
1$ and thus (\ref{nec}) and (\ref{coro}) are satisfied for any $m$ provided
that $p$ is close enough to $1$. Let us mention that this must occur since, at
least when $m^{-}\in L^{\infty}\left(  \Omega\right)  $, $f\left(  \xi\right)
=\xi^{p}$, and either $N=1$ or $N>1$ and $m$ is radial with $0\not \equiv
m\geq0$ in some $B_{r}\left(  0\right)  $, it is known that (\ref{prob}) has a
solution\textit{ }if $p$ is sufficiently close to $1$ (cf. \cite{nodea},
Theorems 2.1 (i) and 3.2). $\blacksquare$

\end{document}